\documentclass{article}
\usepackage{amssymb,amsmath}
 
\input xy 
\xyoption{all}
 
\begin{document}
\title{Hodge decomposition for  elliptic complexes over unital    $C^*$-algebras}
\author{Svatopluk Kr\'ysl \footnote{{\it E-mail address:} Svatopluk.Krysl@mff.cuni.cz, {\it Tel./Fax:} $+$ 420  222 323 221/ $+$ 420 222 323 394} \\{\it \small 
Charles University in Prague, Faculty of Mathematics and Physics
}}
\maketitle \noindent
\centerline{\large\bf Abstract}

For a certain class of complexes of pre-Hilbert $A$-modules, we prove that their
cohomology groups equipped with a canonical quotient structure are again pre-Hilbert $A$-modules and  derive the Hodge decomposition for them. 
We call these complexes self-adjoint parametrix possessing.
We show that $A$-elliptic complexes of differential operator acting on sections of finitely generated projective $A$-Hilbert 
bundles over compact manifolds have this property if the images of certain extensions of their Laplacians are closed.

\section{Introduction}

 In this paper, we focus our attention to specific co-chain complexes of pre-Hilbert modules over unital $C^*$-algebras and adjointable pre-Hilbert
module homomorphisms (differentials) acting between them.  We would like to describe certain properties of their cohomology groups and prove the Hodge decomposition for them.

 Let us recall that a co-chain complex in a category $\mathcal{C}$ is a sequence $d^{\bullet}=(C^k,d_k)_{k\in \mathbb{Z}}$ such
that $C^k$ are  objects in the category $\mathcal{C}$ and $d_k:C^k \to C^{k+1}$ are morphisms in $\mathcal{C}$ which satisfy 
$d_{k+1}  d_k =0,$ $k\in \mathbb{Z}.$  For convenience, we consider complexes bounded
from below, i.e., $k\geq 0.$ As mentioned above, we are interested in complexes in the category of pre-Hilbert $A$-modules and 
pre-Hilbert $A$-module homomorphisms. By a Hodge decomposition, we mean a decomposition of the pre-Hilbert modules $C^k$ in the
given complex into a direct sum of three
pre-Hilbert modules, namely of the module of harmonic elements, the module of closed, and  the module of   co-closed elements. Notice that in particular, 
in order to define the harmonic elements, we have to suppose that the maps $d_k$ forming the complex are adjointable.
Since the cohomology groups of a complex are quotients of the kernel of a map in the complex by the image of the preceding map,
it is not surprising that the cohomology groups does not necessarily belong to the category we started with. In the category of pre-Hilbert modules, the
cohomologies need not be Hausdorff, let alone normed spaces.

  We shall confine ourselves to self-adjoint endomorphisms of  pre-Hilbert modules only. Let $L:V \to V$ be such an endomorphism.
We show that the existence of maps $g, p: V \to V$ such that $1_{|V} = gL + p = Lg + p,$ $Lp=0$ and $p=p^*$  is sufficient
for the decomposition $V= \mbox{Ker}\, L \oplus \mbox{Im}\, L^*$ to hold (no completion involved). 
We call such endomorphisms   self-adjoint parametrix possessing.
Then we apply this result to complexes. 
With each complex $d^{\bullet}=(C^k, d_k)_{k\in \mathbb{N}_0}$ in the category of pre-Hilbert modules and adjointable pre-Hilbert homomorphisms,
 we associate the sequence of self-adjoint endomorphisms
$L_i = d_{i-1} d_{i-1}^* + d_i^* d_i: C^i \to C^i,$ $i\in \mathbb{N}_0,$  the so-called {\it Laplacians} of the  complex. When these Laplacians are self-adjoint parametrix possessing,
we  show that the cohomology groups of the complex are pre-Hilbert modules  isomorphic to $\mbox{Ker}\, L_i,$ and moreover that
$C^i = \mbox{Ker} \, L_i \oplus \mbox{Im} \, d_{i}^* \oplus \mbox{Im} \, d_{i-1}$ (the Hodge decomposition).
 Let us recall that the Hodge theory, i.e., the Hodge decomposition and a description of the cohomology groups, is well
known for complexes of finite dimensional vector spaces and linear maps.
It is a consequence of elliptic operator theory, that the Hodge theory is valid  also for {\it elliptic complexes} of 
differential operators acting between smooth sections of finite rank
real or complex vector bundles over compact manifolds. 

Here, we prove the Hodge decomposition 
for certain kind of complexes $D^{\bullet}=(\Gamma(\mathcal{F}^k), D_k)_{k\in \mathbb{N}_0}$ of $A$-elliptic 
differential operators acting on the smooth sections  of the $A$-Hilbert bundle $\mathcal{F}^i$
(notions introduced by the authors of \cite{FM} and further worked out, e.g., in \cite{ST}).
We  shall suppose  that the operators act between smooth sections of 
Hilbert bundles over compact manifolds, the fibers of the bundles are finitely generated projective Hilbert $A$-modules, 
and that the following technical condition is satisfied:
for each $i\in \mathbb{N}_0,$ the image $\mbox{Im} \, (L_i)_{r_i}$ of the extension $(L_i)_{r_i}:W^{r_i}(\mathcal{F}^i) \to W^{0}(\mathcal{F}^i)$
  of the  Laplacian $L_i$ associated with $D^{\bullet}$
 to the $r_i$-th Sobolev completion  $W^{r_i}(\mathcal{F}^i)$
(of the space of smooth sections of the bundle $\mathcal{F}^i$)  is closed in $W^0(\mathcal{F}^i).$
Here,  $r_i$ denotes the order of the Laplacian $L_i$ and $W^{r_i}(\mathcal{F}^i)$ is the $r_i$-th Sobolev type completion of the space
of smooth sections $\Gamma(\mathcal{F}^i)$ of  $\mathcal{F}^i.$
We derive this application using  \cite{KryslHodge} which is based  on   results of Fomenko, Mishchenko in \cite{FM}.
The authors of the latter article  construct smoothing parametrices for extensions of $A$-elliptic operators to the Sobolev completions. 
In \cite{KryslHodge},   this construction  is used to derive smooth parametrices for $A$-elliptic operators, and a generalization of these results
 to the case of $A$-elliptic complexes. In this sense, the present paper might be considered as a continuation of the work started in \cite{KryslHodge}.

Let us notice that for deriving   results of this paper, we were motivated by the mathematics connected with 
 Quantum theory, in particular,  by results of  Kostant  in \cite{Kostant} and the work of Habermann on the so-called symplectic Dirac operator.
For it, see Habermann, Habermann \cite{H}.
There are also  generalization of the classical Hodge theory (for elliptic operators)
in  directions different from the one we present here. See, e.g., 
Smale et al. \cite{Smale} and the reference there.
For more $K$-theoretically and/or analytically oriented  works, see Pavlov \cite{Pavlov}, Schick
 \cite{Schick}, Shubin \cite{Shubin}, Troitsky \cite{Troitsky}, and Troitsky and Frank \cite{Frank}.

In the second section, we set the terminology and the notation and derive some simple properties of projections, complementability,
and pre-Hilbert module structures on quotients. 
Then we prove that a self-adjoint parametrix possessing endomorphism $L$ of a pre-Hilbert $A$-module $V$ admits a decomposition of the  form
$V = \mbox{Ker}\, L \oplus \mbox{Im} \, L$  (Theorem 3).
In the third section, we derive the Hodge decomposition for self-adjoint parametrix possessing complexes (Theorem 5)
and give a characterization of their cohomology groups (Corollary 7).
In the fourth section, the definitions of an $A$-Hilbert bundle  and an $A$-elliptic complexes are recalled. At the end, 
Theorem 8 on the Hodge theory for the   class of $A$-elliptic complexes mentioned above is proved.

{\bf Preamble:} All manifolds and bundle structures (total spaces, base spaces, and bundle projections) are smooth. 
Base spaces of all bundles are finite dimensional.
Further, if an index of an labeled object exceeds its allowed range, we consider it to be zero.

\section{Parametrix possessing endomorphisms of pre-Hilbert $A$-modules}

 Let $A$ be a unital $C^*$-algebra. We denote the involution, the norm in $A,$ the partial ordering on the double cone of hermitian elements in $A,$ and the unit 
by $^*,$  $| \, |_A,$  $\leq,$ and $1,$ respectively. 

Let us recall that a {\it pre-Hilbert $A$-module} is firstly, a complex vector space $U$ on which $A$ acts.
For definiteness, we consider that $A$ acts from the left, and denote the action by a dot. Secondly, $U$ has to be equipped with a map
$(,)_U: U \times U \to A$ such that for all $a\in A$ and $u,v \in U,$ the following properties hold
\begin{itemize}
\item[1)] $(a.u,v)_U = a^*(u,v)_U$
\item[2)] $(u,v)_U = (v,u)_U^*$
\item[3)] $(u,u)_U \geq 0$ 
\item[4)] $(u,u)_U=0$ if and only if $u=0.$
\end{itemize}
We call such a  map $(,)_U: U \times U \to A$ an $A$-product.
For a pre-Hilbert $A$-module  $(U,(,)_U),$  one defines a norm $| \, |_U: U \to \mathbb{R}^+_0$ (induced by $(,)_U$) by the formula
$U \ni u \mapsto |u|_U = \sqrt{|(u,u)_U|_A}  \in \mathbb{R}_0^+.$
A {\it homomorphism} $L$ between pre-Hilbert $A$-modules $U,V$ has to be  
$A$-linear, i.e., $L(a.u)=a.L(u)$ for each $a\in A$ and 
$u \in U,$ and  continuous  with respect to the norms $| \, |_U$ and $| \, |_V.$
An adjoint of a pre-Hilbert $A$-module homomorphism $L: U \to V$ is a map from $V$ to $U$ denoted by $L^*$ such that
for each $u \in U$ and $v\in V,$  the identity $(Lu, v)_V= (u, L^*v)_U$  holds. If the adjoint exists, it is unique and a pre-Hilbert $A$-module homomorphism.
(See, e.g., Lance \cite{Lance}.) We denote the set of pre-Hilbert $A$-module homomorphisms from $U$ to $V$ by $\mbox{Hom}_A(U,V).$  If $U=V,$ 
$\mbox{End}_A(U)$ denotes $\mbox{Hom}_A(U,V).$
Quite often, in the literature a homomorphism $L: U \to V$ of pre-Hilbert or Hilbert $A$-modules is 
supposed to be adjointable. For technical reasons, we don't follow this convention.
Let us recall, that a pre-Hilbert $A$-module $(U,(,)_U)$ is called a {\it Hilbert $A$-module} if it is complete with respect to $| \, |_U.$

Elements $u,v \in U$ are called orthogonal if $(u,v)_U=0.$ 
For any pre-Hilbert $A$-submodule $U$ of $V,$ we denote by $U^{\bot}$ the orthogonal complement of $U$
defined by $U^{\perp}= \{v\in V \, | (v,u)_V=0 \mbox{ for all } u \in U \}.$
We call $U$ {\it orthogonally complementable} if there exists a pre-Hilbert $A$-submodule $U' \subseteq V$  such that
$U \oplus U' = V.$ 
Let us notice that if we write a direct sum of pre-Hilbert $A$-submodules, we suppose that the elements belonging to different
summands are mutually orthogonal. 
It is immediate to see that for any pre-Hilbert $A$-submodules $V\subseteq W$ of a pre-Hilbert 
$A$-module $U,$ the operation of taking the orthogonal complement changes the inclusion sign, i.e., 
\begin{eqnarray}
V^{\bot} \supseteq W^{\bot} \label{anti}
\end{eqnarray}

\subsection{Complementability, quotients and parametrix possessing maps}

Let us start with the following  simple observation.
For any pre-Hilbert $A$-module $V,$   
an element $p$ of $\mbox{End}_A(V)$ is called a projection if
$p^2 = p.$
 Let  $p$ be a projection and let us denote by $U$ the $A$-submodule $\mbox{Im}\, p.$ 
For each $z \in U,$ there exists $x\in V$ such that $px=z.$ Thus, 
$p^2x=pz$ which in turn implies $pz=px=z,$ i.e., if $p$ is a projection onto an
$A$-submodule $U,$ then if restricted to $U,$ $p$ is
the identity on $U.$ If $V=U\oplus U'$ and $p(x_U+ x_{U'})=x_U,$ where $x_U \in U$ and $x_{U'}\in U',$ we call $p$ a projection onto
$U$ along $U'.$

{\bf Lemma 1:} Let $V$ be a pre-Hilbert $A$-module and  $U$ be an orthogonally complementable pre-Hilbert $A$-submodule of $V.$ 
If $V= U \oplus U',$ then $U'=U^{\bot}$
and the projection $p$ onto $U$ along $U^{\bot}$ is self-adjoint.
Conversely, if $p$ is a self-adjoint projection in $V,$ then $U=\mbox{Im} \, p$ is an orthogonally complementable pre-Hilbert $A$-submodule 
and $1-p$ is a projection onto $U^{\bot}$ along $U.$

{\it Proof.}
For $x \in U^{\bot},$ there are uniquely determined $x_U\in U$ and $x_{U'} \in U'$ for which $x= x_U + x_{U'}.$
Let us compute $(x_U,x_U)_V = (x-x_{U'},x_U)_V=(x,x_U)_V -(x_{U'},x_U)_V= (x,x_U)_V=0$ since $x\, \bot \, U.$ Thus, $x_U=0$ proving
$U^{\bot}\subseteq U'.$
 The opposite inclusion follows from the definition of the orthogonal complement immediately.
Further, for any $x\in V$ and $y=y_U+y_{U'}\in V, y_U\in U, y_{U'}\in U',$ let us write
$(px,y)_V = (x_U, y_U+y_{U'})_V= (x_U, y_U)_V = (x, y_U)_V=(x, py)_V,$ i.e., $p$ is self-adjoint.  

For the other statement, set $U=p(V)$ and $U'=(1-p)(V).$ We have $V=U + U'.$ 
For $x \in U \cap U',$  we get $(x,x)_V=(px,(1-p)x)_V = (px,x)_V - (px,px)_V= (x,x)_V-(x,x)_V=0,$ i.e., $x=0$ and thus, the sum is direct.
It is immediate to see that $1-p$ is self-adjoint and a projection.
 To prove that $1-p$ projects onto $U^{\bot},$ let  us consider an element
$y \in U^{\bot}.$ We may compute $(py,py)_V = (y,p^*py)_V = (y,p^2y)_V = (y,py)_V= 0.$ Thus, $py = 0$ and therefore
$y = (1-p)y.$ Obviously, $1-p$ annihilates the elements from $U.$
Thus, $1-p$ is a self-adjoint projection onto $U^{\bot}$ along $U.$
$\Box$

For a normed space $(Y,\, |\,|_Y)$ and its closed normed subspace
$X$, one usually considers  the quotient space $Y/X$ equipped with the norm $| \, |_q: Y/X \to \mathbb{R}$ 
defined by $$|[y]|_q = \mbox{inf}\{|y-x|_Y, x \in X \},$$  where $y\in Y$
and $[y]$ denotes the equivalence class of $y.$
We call $|\,|_q$ the quotient norm.
It is immediate to see that if $Y$ is a Banach space, the quotient is a Banach space as well.
Now, we focus our attention to quotients of pre-Hilbert $A$-modules.
When we speak of a quotient $V/U$ of a pre-Hilbert module $V$ and  its orthogonally complementable submodule $U$, we 
think  of $V/U$ as of an $A$-module equipped with the following  $A$-product $(,)_{V/U}.$
Let $p$ be the projection onto $U^{\bot}$ along $U.$ We set
$([u],[v])_{V/U} =(p(u),p(v))_V,$  $u,v \in V.$
The map $(,)_{V/U}$ is easily seen to be correctly defined. Further, it is evident that it maps into the set of non-negative elements of $A.$
Suppose that $([u],[u])_{V/U}=0$ for an element $u\in V.$ Then $(p(u),p(u))_V = 0$ and consequently, $p(u)=0.$ Thus $u\in U,$ i.e.,
$[u]=0,$ proving the positive definiteness of  $(,)_{V/U}$ (properties 3 an 4 in the definition of the  $A$-product).
 Summing up, in the case of an orthogonally complementable pre-Hilbert $A$-submodule $U$  of a pre-Hilbert 
$A$-module $V,$ we obtain a canonical pre-Hilbert $A$-module structure $(V/U, (,)_{V/U}).$

{\bf Lemma 2:} Let $U$ be an orthogonally complementable pre-Hilbert $A$-sub\-mo\-dule of a pre-Hilbert $A$-module $(V, (,)_V).$
Then  
\begin{itemize}
\item[1)] $V/U$ and $U^{\bot}$ are isomorphic pre-Hilbert $A$-modules and 
\item[2)] the quotient norm $| \, |_q$ coincides with the norm induced by $(,)_{V/U}.$
\end{itemize}

{\it Proof.}  Let $p$ be the projection onto $U^{\bot}$  along $U$ and $p' =1-p$ be the projection onto
$U$ along $U^{\bot}$ (Lemma 1).
For any $v\in V,$ we have 
\begin{eqnarray*}
|[v]|_q^2 &=& \mbox{inf}_{u\in U}|v-u|_V^2 \\
&=& \mbox{inf}_{u \in U}|(v-u,v-u)_V|_A\\
&=& \mbox{inf}_{u\in U}|(p'v + pv-u,p'v +pv -u)_V|_A \\
&=& \mbox{inf}_{u\in U}|(p'v - u, p'v +pv - u)_V + (pv, p'v + pv -  u)_V|_A \\
&=& \mbox{inf}_{u\in U} |(p'v - u, p'v -u)_V + (pv,pv)_V|_A\\
&=& |(pv,pv)_V|_A = |[v]|_{V/U}^2,
\end{eqnarray*}
where in the second last step, we used the fact that $|a+b|_A\geq |a|_A$ which holds for any non-negative $a,b \in A.$
This proves the second item.

It is easy to check that $\Phi([v])=pv$ is a well defined $A$-module homomorphism of $V/U$ into $U^{\bot}.$ 
Consider also the map $\Psi: U^{\bot} \to V/U$ defined by $\Psi(u)=[u],$ $u \in U^{\bot}.$  Both of the maps are
continuous with respect to the norm topology on $U^{\bot}$ (inherited from $(V, | \, |_V)$) and the quotient topology on $V/U.$
Because  the topology induced by $|\,|_q$ coincides  with the quotient topology, and $|\,|_q$ coincides with $|\,|_{V/U}$ (due to the first paragraph of this proof),
 we conclude that both $\Phi$ and $\Psi$ are continuous
with respect to the norms $| \, |_{U^{\bot}}$ and $|\,|_{V/U},$ i.e., they are homomorphisms of the corresponding
pre-Hilbert $A$-modules. Here, $|\,|_{U^{\bot}}$ denotes the restriction of $|\,|_V$ to $U^{\bot}.$
Further, for any $u\in U^{\bot},$ we have $\Phi(\Psi(u))=\Phi([u])=pu=u$ since $p$ projects onto $U^{\bot}.$
For each $[v]\in V/U,$ we may write $\Psi(\Phi([v]))=\Psi(pv)=[pv].$ Because the difference of $v$ and $pv$ lies in $U,$
we get  $\Psi \circ \Phi = 1_{|V/U}.$ 
$\Box$

{\bf Remark 1:} As a consequence of Lemma 2, for a pre-Hilbert $A$-module $V$ and an orthogonally complementable 
pre-Hilbert $A$-submodule $U$ of $V$ if  $(V/U, | \, |_q)$ is a Banach space, $(V/U, (,)_{V/U})$ is a Hilbert 
$A$-module. Further if moreover, $V$ is Hilbert $A$-module, then $(V/U, (,)_{V/U})$ is a Hilbert $A$-module as well.

Now, we shall  focus our attention to  relations  between orthogonal complementability of images of pre-Hilbert $A$-module endomorphisms
and the property described in the next definition.


{\bf Definition 1:} Let $L$ be an endomorphism of a pre-Hilbert module $V.$ We call $L$
{\it parametrix possessing} if there exists pre-Hilbert $A$-module endomorphisms $p,g: V \to V$ such that
\begin{eqnarray*}
1_{|V} &=& g L + p\\
1_{|V} &=& Lg + p\\
L p &=& 0
\end{eqnarray*}
where $1_{|V}$ denotes the identity on $V.$
We call a parametrix possessing map $L$ {\it self-adjoint parametrix possessing} if $L$ and $p$ are self-adjoint.

{\bf Remark 2:} The first two equations in Definition 1 are  called parametrix equations.
Notice that there exist pre-Hilbert $A$-module endomorphisms which are not parametrix possessing and also such for which,
$g$ and $p$ are not uniquely determined. The name parametrix is borrowed from the theory of elliptic PDEs.

{\bf Theorem 3:} Let $L: V \to V$ be a  self-adjoint  parametrix possessing endomorphism of a pre-Hilbert 
$A$-module $V$ with the corresponding maps denoted by $g$ and $p.$  Then 
\begin{itemize}
\item[1)] $p$ is a projection onto $\mbox{Ker} \, L$ and
\item[2)] $V = \mbox{Ker} \,  L \oplus \mbox{Im} \,  L.$
\end{itemize}

{\it Proof.} 
\begin{itemize}
\item[1)]
Composing the first parametrix equation from the right by $p$ and using the third equation
from the definition of a parametrix possessing endomorphism, we get that $p^2=p,$ i.e., $p$ is an
idempotent. 
Restricting $1_{|V} = gL + p$ 
to $\mbox{Ker} \,  L,$ we get $1_{|\mbox{Ker} \,  L} = p_{|\mbox{Ker} \,  L}$ which implies that $\mbox{Im} \,  p \supseteq \mbox{Ker} \,  L.$
Further, $L p=0$ forces $\mbox{Im} \,  p \subseteq \mbox{Ker} \,  L.$ Thus, $\mbox{Im}\, p =\mbox{Ker} \, L$ and consequently, $p$ is a projection onto $\mbox{Ker} \,  L.$
\item[2)]
Since
$p$ is self-adjoint, we may use Lemma 1 to conclude that
$V= \mbox{Ker} \,  L \oplus (\mbox{Ker} \,  L)^{\bot}.$ It is sufficient  to prove the equality
\begin{eqnarray}
\mbox{Im} \,  L = (\mbox{Ker} \,  L)^{\bot} \label{bot}
\end{eqnarray}
 First, we prove that $\mbox{Im} \,  L \subseteq (\mbox{Ker} \,  L)^{\bot}.$
Let $y = Lx$ for an element $x \in V.$
For any $z\in \mbox{Ker} \,  L,$ we may write $(y,z) = (Lx,z) = (x,L^*z) = (x,Lz) =  0.$ Thus, 
$y \, \bot \, \mbox{Ker} \,  L.$
Now, we prove that $(\mbox{Ker} \,  L)^{\bot} \subseteq \mbox{Im} \,  L.$ Let $x \in (\mbox{Ker} \,  L)^{\bot}.$
Using the second parametrix equation, we obtain $Lg x = (1-p)x = x$ since $1-p$ projects onto
$(\mbox{Ker} \,  L)^{\bot}$  (Lemma 1).
Therefore $x = Lgx \in \mbox{Im} \,  L.$
Summing up, $\mbox{Im} \,  L = (\mbox{Ker} \,  L)^{\bot},$ and the equation $V = \mbox{Ker} \,  L \oplus \mbox{Im} \,  L$ follows.
\end{itemize}
$\Box$

{\bf Remark 3:} Let us notice that in particular, any 
self-adjoint parametrix possessing endomorphism has  closed image.

\section{Cohomology and Hodge decomposition}

In this section, we focus our attention to co-chain complexes $d^{\bullet} = (C^k, d_k)_{k\in \mathbb{N}_0}$ of pre-Hilbert $A$-modules
and adjointable pre-Hilbert $A$-module homomorphisms, i.e., for each $k\in \mathbb{N}_0,$ $d_k: C^k \to C^{k+1}$ is an adjointable pre-Hilbert $A$-module homomorphism 
which satisfies $d_{k+1} d_k=0.$
We will  transfer Theorem 3 to the situation of co-chain complexes.
Let us consider
the sequence of {\it Laplacians} $L_k=d_k^*d_k + d_{k-1}d_{k-1}^*,$ $k\in \mathbb{N}_0,$ associated with $d^{\bullet}.$

{\bf Definition 2:} Let $d^{\bullet}=(C^{k}, d_{k})_{k \in \mathbb{N}_0}$ be a co-chain complex of pre-Hilbert $A$-modules
  and adjointable pre-Hilbert $A$-module homomorphisms.
  We call  $d^{\bullet}$ a {\it parametrix possessing complex} if for each $k\in \mathbb{N}_0,$ the  
associated Laplacian $L_k$ is a parametrix possessing endomorphism of $C^k$. 
We call $d^{\bullet}$ a self-adjoint parametrix possessing complex if 
the operators $L_k$ are self-adjoint parametrix possessing pre-Hilbert $A$-module 
endomorphisms for all $k\in \mathbb{N}_0.$ 

Obviously,  a parametrix possessing complex is self-adjoint parametrix 
possessing if and only if the projection $p_k$ is self-adjoint, $k\in \mathbb{N}_0.$
Notice that (in concordance with the preamble), $L_0=d_0^*d_0.$

{\bf Lemma 4:}  Let $d^{\bullet}=(C^{k}, d_{k})_{k\in \mathbb{N}_0}$ be a co-chain complex of pre-Hilbert $A$-modules and
adjointable pre-Hilbert $A$-module homomorphisms. Then
$$\mbox{Ker} \, L_k = \mbox{Ker} \, d_k \cap \mbox{Ker} \, d_{k-1}^*$$

{\it Proof.} The inclusion $\mbox{Ker}\, L_k \supseteq \mbox{Ker} \, d_k \cap \mbox{Ker} \, d_{k-1}^*$  follows from the definition
of the Laplacian $L_k.$ To prove the opposite one,
let $x \in \mbox{Ker} \, L_k.$ We may write  
$0 =(x,L_k x)_{C^k} = (x, d_k^*d_k x + d_{k-1}d_{k-1}^*x)_{C^k} 
  = (d_k x, d_k x)_{C^{k+1}} + (d_{k-1}^*x,d_{k-1}^*x)_{C^{k-1}}.$  
Thus, $d_k x = d_{k-1}^* x =0$ due to the positive definiteness of the $A$-products on $C^{k+1}$ and $C^{k-1},$ respectively. $\Box$




{\bf Theorem 5:}
Let $d^{\bullet} = (C^{k}, d_{k})_{k\in \mathbb{N}_0}$ be a self-adjoint parametrix possessing complex.
Then for any $k\in \mathbb{N}_0,$ we have the decomposition
$$C^k = \mbox{Ker} \, L_k \oplus \mbox{Im} \,  d_{k-1} \oplus \mbox{Im} \,  d_k^*.$$

{\it Proof.} Because  $d^{\bullet}$ is a  parametrix possessing complex, there exist maps
$g_k$ and $p_k$ satisfying the parametrix equations for $L_k$ and the identity $L_k p_k = 0,$   $k\in \mathbb{N}_0.$

\begin{itemize}
\item[1)]
Due to Lemma 4, we have $\mbox{Ker} \, L_k \subseteq \mbox{Ker}\, d_{k-1}^*.$ Therefore using the formulas (\ref{anti}) and (\ref{bot}), 
we get  $(\mbox{Ker} \,  d_{k-1}^*)^{\bot} \subseteq (\mbox{Ker} \,  L_k)^{\bot}  = \mbox{Im} \, L_k.$
Further, due to Lemma 4 again, we have  $\mbox{Ker} \,  L_k \subseteq  \mbox{Ker}\,d_k.$ In the same way as above, we get
$(\mbox{Ker} \,  d_k)^{\bot} \subseteq (\mbox{Ker} \,  L_k)^{\bot} = \mbox{Im} \,  L_k.$ 
Summing up, $(\mbox{Ker} \,  d_{k-1}^*)^{\bot} + (\mbox{Ker} \,  d_k)^{\bot}   \subseteq \mbox{Im} \, L_k.$ 

\item[2)]
The inclusion  $\mbox{Im} \,  d_{k-1} \subseteq (\mbox{Ker} \,  d_{k-1}^*)^{\bot}$ holds since for
any $x\in C^{k-1}$ and $y \in \mbox{Ker} \,  d_{k-1}^*,$ we have
$(d_{k-1} x, y)_{C^k} = (x,d_{k-1}^*y)_{C^{k-1}} = 0.$ Similarly, 
$\mbox{Im} \,  d_{k}^* \subseteq  (\mbox{Ker} \,  d_k)^{\bot}.$ 
Combining this with the result of item 1 of this proof, we get  $\mbox{Im} \, d_{k-1} + \mbox{Im} \, d_k^* \subseteq (\mbox{Ker} \,  d_{k-1}^*)^{\bot} 
+ (\mbox{Ker} \,  d_k)^{\bot}  \subseteq \mbox{Im} \,  L_k.$
Now, we show that the sum $\mbox{Im} \, d_k^* + \mbox{Im} \, d_{k-1}$ is direct. Let $y = d_k^*x = d_{k-1}z$ for elements $x \in C^{k+1}$ 
and $z \in C^{k-1}.$
We have $(y,y)_{C^k} = (d_k^*x, d_{k-1}z)_{C^k} = (x,d_k d_{k-1}x)_{C^{k+1}} = 0,$ and consequently, $y=0.$ Summing up,
$\mbox{Im} \, d_k^* \oplus \mbox{Im} \, d_{k-1} \subseteq \mbox{Im}\, L_k.$

\item[3)]
It is easy to prove that $\mbox{Im}\, L_k \subseteq \mbox{Im}\, d_k^* \oplus \mbox{Im} \, d_{k-1}.$ 
Indeed, for any $y\in \mbox{Im} \,  L_k,$ there exists $x\in C^k$ 
such that  $y = L_k x = d^*_k d_k x + d_{k-1} d^*_{k-1} x = d_k^* (d_kx) + d_{k-1}(d_{k-1}^*x) 
\in \mbox{Im} \, d^*_k + \mbox{Im} \,  d_{k-1}.$ This together with item 2 proves that
$\mbox{Im}\, L_k =  \mbox{Im}\, d_k^* \oplus \mbox{Im} \, d_{k-1}.$

\item[4)] Because $L_k$ is a self-adjoint parametrix possessing pre-Hilbert $A$-module 
endomorphism of $C^k,$ we have
due to Theorem 3, the equality  $C^k = \mbox{Im} \, L_k \oplus \mbox{Ker} \, L_k.$ Substituting for
$\mbox{Im} \, L_k$ from item 3 of this proof, we get the  decomposition.
\end{itemize}
$\Box$

{\bf Remark 4:} During the proof of the previous theorem, we obtained for (a
self-adjoint parametrix-possessing complex) $d^{\bullet},$ the decomposition
$$\mbox{Im} \, L_k =  \mbox{Im}\, d_k^* \oplus \mbox{Im} \, d_{k-1}.$$

 Notice that if $d^{\bullet}=(C^k, d_k)_{k\in \mathbb{N}_0}$ is a co-chain complex, then its adjoint
$(C^k, d_k^*)_{k\in \mathbb{N}_0}$ is a chain complex 
as one easily sees from $d_{k}^*d_{k+1}^* = (d_{k+1}d_k)^*.$

{\bf Theorem  6:} Let $d^{\bullet} = (C^{k}, d_{k})_{k\in \mathbb{N}_0}$ be a self-adjoint parametrix possessing complex.
Then for any $k\in \mathbb{N}_0,$ 
\begin{eqnarray*}
\mbox{Ker} \, d_k &=& \mbox{Ker} \, L_k \oplus \mbox{Im} \, d_{k-1}\\
\mbox{Ker} \, d_k^* &=& \mbox{Ker} \, L_{k+1} \oplus \mbox{Im} \, d_{k+1}^*.
\end{eqnarray*}

{\it Proof.} 
Because of Theorem 5, we know that the sums in both rows are direct.

The inclusion $\mbox{Ker} \, L_k \oplus \mbox{Im} \, d_{k-1} \subseteq \mbox{Ker} \, d_k$ is an immediate consequence of the definition 
of a co-chain complex and of Lemma 4. To prove the opposite inclusion, let us consider an element $y \in \mbox{Ker} \, d_k.$ 
Due to Theorem 5, there exist elements $y_1 \in \mbox{Ker} \, L_k,$ $y_2 \in \mbox{Im} \, d_{k-1}$ and $y_3 \in \mbox{Im} \, d_{k}^*$ such that
$y = y_1 + y_2 + y_3.$ It is sufficient to prove that $y_3 = 0.$ Let $z_3 \in C^{k+1}$ be such that $y_3 =d_k^* z_3.$
We have $0 = (d_k y, z_3) = (d_k y_1 + d_k y_2 + d_k y_3, z_3) = (d_k y_3, z_3)= (y_3, d_k^* z_3) = (y_3, y_3)$ which implies 
$y_3 = 0,$ and the first relation follows.

The inclusion $\mbox{Ker} \, L_{k+1} \oplus \mbox{Im} \, d_{k+1}^* \subseteq \mbox{Ker} \, d_k^*$ follows from Lemma 4 and the second part of Remark 4.
To prove the inclusion $\mbox{Ker} \, d_k^* \subseteq \mbox{Ker} \, L_{k+1} \oplus \mbox{Im} \, d_{k+1}^*,$  we proceed similarly as in the previous paragraph.
For $y\in \mbox{Ker} \, d_k^*,$ we have $y_1 \in \mbox{Ker}\, L_{k+1},$ $y_2 \in \mbox{Im}\, d_k,$ and $y_3 \in \mbox{Im} \, d_{k+1}^*$ such that
$y=y_1+y_2+y_3$ (Theorem 5). Let us consider an element $z_2 \in C^k$ for which $y_2=d_k z_2.$
We have $0=(d_k^*y,z_2)=(d_k^*y_1+d_k^*y_2+d_k^*y_3, z_2) = (d_k^*y_2,z_2)=(y_2,y_2).$ Thus,
$y_2 = 0,$ and the second relation follows.
$\Box$

For any complex $d^{\bullet} = (C^k, d_k)_{k\in \mathbb{N}_0},$ we consider the cohomology groups
$$H^i(d^{\bullet}, A) = \frac{\mbox{Ker} \, (d_i: C^i \to C^{i+1})}{\mbox{Im} \, (d_{i-1}: C^{i-1} \to C^i)}, i \in \mathbb{N}_0.$$
Notice that in general, the $A$-module $Z^i(d^{\bullet},A)=\mbox{Im} \, (d_{i-1}:C^{i-1} \to C^i)$ of co-boundaries needs not be
orthogonally complementable or even not closed in the appropriate pre-Hilbert $A$-module. Thus, the cohomology group needn't be a Hausdorff space with respect 
to the quotient topology.
Nevertheless, in the case of a self-adjoint parametrix possessing complex, we derive

{\bf Corollary 7:} Let $d^{\bullet} = (C^{k}, d_k)_{k\in \mathbb{N}_0}$ be a self-adjoint parametrix possessing complex.
Then for each $i,$ the cohomology group  $H^i(d^{\bullet},A)$ is a pre-Hilbert $A$-module.
If $d^{\bullet}$ is a self-adjoint parametrix possessing complex of Hilbert $A$-modules, then 
for each $i,$ the cohomology group  $H^i(d^{\bullet},A)$ is a Hilbert $A$-module.

{\it Proof.}  Because of Theorem 6, $U=\mbox{Im} \, d_{i-1}$ is orthogonally complementable in $V=\mbox{Ker}\, d_i.$ 
Thus using Lemma 2, the cohomology group $H^i(d^{\bullet}, A) = \mbox{Ker} \, d_i / \mbox{Im} \, d_{i-1}$ 
equipped with the canonical $A$-product $(,)_{V/U}$ is a pre-Hilbert $A$-module.
The second statement follows from  Remark 1. $\Box$

{\bf Remark 5:} Notice that moreover due to Theorem 6, we have
$$H^i(d^{\bullet}, A) \cong \mbox{Ker} \, L_i$$ as pre-Hilbert $A$-modules.

\section{Application to differential operators}

Let $M$ be a   finite dimensional manifold and $p:\mathcal{F} \to M$ be a Banach bundle over $M.$
We call $p:\mathcal{F} \to M$ an {\it $A$-Hilbert bundle} if there exists a Hilbert $A$-module  $(S, (,)_S)$ and a bundle atlas $\mathcal{A}$ of $p$ compatible
with the bundle atlas of $p$ considered as the Banach bundle only, such that
\begin{itemize}
\item[1)] $(S, | \, |_S)$ is the typical fiber of the Banach bundle $p$
\item[2)] for each $m\in M,$ the fiber $\mathcal{F}_m=p^{-1}(m)$ is equipped with a Hilbert $A$-product, denoted by $(,)_m$ 
\item[3)] for each $m\in M$ and each chart $(\phi_U,U) \in  \mathcal{A},$ $U\ni m,$ the map
${\phi_U}_{|\mathcal{F}_m}: (\mathcal{F}_m, (,)_m) \to (S,(,)_S)$ is a Hilbert $A$-module isomorphism and 
\item[4)] the transition maps of the charts in $\mathcal{A}$ are maps into $\mbox{Aut}_A(S).$
\end{itemize}
Let us recall that for two bundle charts $\phi_U: p^{-1}(U) \to U \times S$ and $\phi_V: p^{-1}(V) \to V \times S,$ 
their transition map $\phi_{UV}:U \cap V \to \mbox{Aut}_A(S)$ is  defined by the formula
$(\phi_U \circ \phi_V^{-1})(m, v) = (m, \phi_{UV}(m)v)$ for each $m \in U \cap V$ and $v \in S.$
A {\it homomorphism of $A$-Hilbert bundles} $p_1:\mathcal{F}_1 \to M$ and $p_2: \mathcal{F}_2 \to M$ is a map
$F: \mathcal{F}_1 \to \mathcal{F}_2$ between the total spaces of 
$p_1$ and $p_2$ such that $p_2 \circ F = p_1,$ and in each fiber, $F$ is a Hilbert $A$-module homomorphism, i.e., for each $m\in M,$ $F_{|p_1^{-1}(m)}: 
(\mathcal{F}_1)_m \to (\mathcal{F}_2)_m$ is such a map. 
An $A$-Hilbert bundle is called finitely generated projective
if and only if its typical fiber, the Hilbert $A$-module $(S,(,)_S),$ is a finitely generated  projective Hilbert $A$-module.

Let us suppose that $M$ is compact, equipped with a Riemannian metric $g$ and let us choose a volume element $|\mbox{vol}_g|$ of $(M,g).$ 
The (positive definite) Laplace-Beltrami operator will be denoted by $\triangle_g.$
 For each $A$-Hilbert bundle $p:\mathcal{F} \to M$ over $M$ and each $t\in \mathbb{Z},$ Fomenko and Mishchenko in \cite{FM} define a certain Hilbert $A$-module,
the so-called Sobolev completion $W^t(\mathcal{F})$ of the space of smooth sections $\Gamma(\mathcal{F})$ of $\mathcal{F}$. Let us sketch their construction 
briefly.
Obviously, the space $\Gamma(\mathcal{F})$ of smooth sections of $\mathcal{F}$ carries a left $A$-module structure given by
$(a.s)(m) = a.(s(m)),$ $a\in A,$ $s\in \Gamma(\mathcal{F})$ and $m\in M.$
One defines an $A$-product by the formula
$$(s,s')_{\Gamma} = \int_{m\in M} (s, s')_m |\mbox{vol}_g|_m, \mbox{  } s'\in \Gamma(\mathcal{F}).$$
Setting $$(s,s')_t = \int_{m\in M} (s,(1+\triangle_g)^ts')_m |\mbox{vol}_g|_m, \mbox{  } s'\in \Gamma(\mathcal{F}),$$
we obtain further pre-Hilbert $A$-modules $(\Gamma(\mathcal{F}),(,)_t).$ Obviously, $(,)_{\Gamma}= (,)_0.$
For definiteness, we consider (the appropriate manifold version of) the Bochner integral of Banach space valued functions.
The Sobolev completion $W^t(\mathcal{F})$ is defined as the completion of $\Gamma(\mathcal{F})$ with respect to the norm $| \, |_t$ induced by
the $A$-product $(,)_t.$ We keep denoting the Hilbert $A$-products by $(,)_t$ also if we consider their extensions to
$W^t(\mathcal{F}).$ See Fomenko, Mishchenko \cite{FM} or Solovyov, Troitsky \cite{ST} for details on this construction if necessary.

Our reference for the statements in the upcoming paragraph is Solovyov, Troitsky \cite{ST}.
For a definition of an $A$-differential operators we refer to Solovyov, Troitsky \cite{ST}, pp. 79 and 80.
We omit the prefix $A$- and call these operators differential operators only.
 For any differential operator $D: \Gamma(\mathcal{F}_1) \to \Gamma(\mathcal{F}_2),$ 
we have the order $\mbox{ord}(D) \in \mathbb{Z}$ of $D$, the adjoint $D^*:\Gamma(\mathcal{F}_2)
 \to \Gamma(\mathcal{F}_1)$ (Theorem 2.1.37 in \cite{FM}), and for each $t\in \mathbb{Z},$ 
the  (continuous) extension $D_t: W^t(\mathcal{F}_1) \to W^{t - \mbox{ord}(D)}(\mathcal{F}_2)$ of $D$ at our disposal
(\cite{ST}, pp. 89, Theorem 2.1.60). 
 Let us denote by $\pi: T^*M \to M$ the  cotangent bundle and let $\pi'$ be  the restriction of $\pi$ to
$T^*M'=T^*M \setminus \{(m,0) \in T^*M | m\in M\}.$
For a differential operator  $D,$ one defines the notion of its symbol $\sigma(D):\pi^*(\mathcal{F}_1) \to \mathcal{F}_2.$
See \cite{ST} pp. 79 and 80.
If $T^*M$ is considered with the trivial $A$-Hilbert bundle structure, i.e., $a . \alpha_m = \alpha_m$ for each $a\in A$ and $\alpha_m \in T^*_mM,$ $m\in M,$
the symbol $\sigma(D): \pi^*(\mathcal{F}_1) \to \mathcal{F}_2$ is an adjointable $A$-Hilbert bundle homomorphism. 
The restriction of the symbol $\sigma=\sigma(D)$ of $D$ to ${\pi'}^*(\mathcal{F}_1)$ will be denoted by
$\sigma'.$ 

Let  $(p_k: \mathcal{F}^k \to M)_{k\in \mathbb{N}_0}$ be 
a sequence of finitely generated projective $A$-Hilbert bundles over $M$ and
$D^{\bullet} = (\Gamma(\mathcal{F}^k),D_k)_{k\in \mathbb{N}_0}$ be  a complex of differential operators in $\mathcal{F}^{\bullet},$
i.e., $D_k:\Gamma(\mathcal{F}^k) \to \Gamma(\mathcal{F}^{k+1})$ is a differential operator and $D_{k+1}  D_k  =0,$ $k\in \mathbb{N}_0.$ Let  us
set  $\sigma_{k} = \sigma(D_k)$ for the symbol of $D_{k}.$
The symbol sequence $\sigma^{\bullet} = (\pi^*(\mathcal{F}^k), \sigma_k)_{k\in \mathbb{N}_0}$ is a 
complex in the category of $A$-Hilbert bundles. 

{\bf Definition 3:} A complex $D^{\bullet} = (\Gamma(\mathcal{F}^k), D_k)_{k \in \mathbb{N}_0}$ 
of differential operators in $A$-Hilbert bundles is called {\it $A$-elliptic} if its restricted symbol sequence 
${\sigma'}^{\bullet} = (\sigma'_k, {\pi'}^*(\mathcal{F}^k))_{k\in \mathbb{N}_0}$ is an   exact sequence  in the category of $A$-Hilbert bundles.

 According to classical conventions, we denote the Laplacians 
$L_k$ associated with a complex $D^{\bullet}$ of differential operators by $\triangle_k.$
Let $r_k$ denote the order of $\triangle_k$.  

{\bf Remark 6:} 
\begin{itemize}
\item[1)] If $D:\Gamma(\mathcal{E}) \to \Gamma(\mathcal{F})$ is a single differential operator, then we consider it
as the complex $$0\to \Gamma(\mathcal{E}) \overset{D}{\rightarrow} \Gamma(\mathcal{F}) \to 0.$$
In this case, the definition of an $A$-elliptic complex of differential operators coincides with the (classical)
definition of an $A$-elliptic operator given in Solovyov, Troitsky \cite{ST}.
\item[2)]If $D^{\bullet}$ is an $A$-elliptic complex of differential operators, then for each $i\in \mathbb{N}_0,$ the Laplacian 
$\triangle_i$ is an $A$-elliptic operator. For it, see Corollary 10 in Kr\'ysl \cite{KryslHodge}.
\end{itemize}

Let us state the following

{\bf Theorem 8:} Let $A$ be a unital $C^*$-algebra and $D^{\bullet} = (\Gamma(\mathcal{F}^k), D_k)_{k \in \mathbb{N}_0}$ be 
an $A$-elliptic complex in finitely generated projective $A$-Hilbert bundles over a compact manifold $M$. Let for each $k\in \mathbb{N}_0,$
 the image of the extension $(\triangle_k)_{r_k}$ of $\triangle_k$ be closed in $W^{0}(\mathcal{F}^k).$ 
Then for any $i\in \mathbb{N}_0,$
\begin{itemize}
\item[1)] $H^i(D^{\bullet},A)$ is a finitely generated projective Hilbert $A$-module
\item[2)] $\Gamma(\mathcal{F}^i) = \mbox{Ker} \, \triangle_i \oplus \mbox{Im} \, D_i \oplus \mbox{Im} \, D_{i-1}^*$
\item[3)] $\mbox{Ker} \, D_i = \mbox{Ker} \, \triangle_{i} \oplus \mbox{Im} D_i^*$
\item[4)] $\mbox{Ker} \, D_i^* = \mbox{Ker} \, \triangle_{i+1} \oplus \mbox{Im} \, D_i $
\end{itemize}

{\it Proof.} 
 In the proof of Theorem 8 in \cite{KryslHodge}, a projection $P$ is constructed which 
satisfies  the parametrix equations for a self-adjoint $A$-elliptic operator $K:\Gamma(\mathcal{F}) \to \Gamma(\mathcal{F})$ of order $r$  provided  $\mbox{Im}\, K_r$ is closed in $W^0(\mathcal{F}).$ 
(Recall that $K_r$ denotes the extension of $K$ to $W^r(\mathcal{F}).$)
Let us sketch this construction briefly.
For $K_{r}: W^{r}(\mathcal{F}) \to W^0(\mathcal{F}),$  consider its adjoint
 $(K_{r})^*:W^0(\mathcal{F}) \to W^{r}(\mathcal{F})$
and the projection $p_{\mbox{Ker} \, K_r^*}$  from the space  $W^0(\mathcal{F})$  onto the kernel
$\mbox{Ker} \, (K_r)^*.$  Let us denote the restriction of  $p_{\mbox{Ker} \, K_r^*}$ to $\Gamma(\mathcal{F})$ by $P.$
Besides $P,$ a further map  $G:\Gamma(\mathcal{F}) \to \Gamma(\mathcal{F})$
is constructed  in the mentioned proof,  which satisfies
$1_{|\Gamma(\mathcal{F})}=GK + P = KG  + P$ and $Kp=0.$ (See Theorem 8 in \cite{KryslHodge}.)  In particular, $K$ is a parametrix possessing  endomorphism of
 the pre-Hilbert $A$-module $(\Gamma(\mathcal{F}), (,)_{\Gamma}).$   

Now, let us prove that the operator $K$ is also self-adjoint parametrix possessing.
Due to the mentioned closed image assumption on the extension $K_r$ of $K,$
$\mbox{Ker} \, (K_r)^*$ is orthogonally complementable in $W^0(\mathcal{F})$  according to Lance \cite{Lance}, Theorem 3.2,  pp. 22.
Thus, the projection $p_{\mbox{Ker} \, K_r^*}$ is self-adjoint according to Lemma 1.
Restricting $p_{\mbox{Ker} \, K_r^*}$ to $\Gamma(\mathcal{F})$ does not change the
property of being idempotent. Moreover, the restriction keep being self-adjoint,  because 
the $A$-product $(,)_{\Gamma}$ in $\Gamma(\mathcal{F})$ coincides with the $A$-product $(,)_0$ in $W^0(\mathcal{F})$  when restricted to $\Gamma(\mathcal{F})\times \Gamma(\mathcal{F}).$
Thus, $K$ is not only parametrix possessing, but it is also self-adjoint parametrix possessing.

Now, let us pass to the statement we are proving. Let $i \in \mathbb{N}_0.$ Using item 2 of Remark 6, $\triangle_i$  is $A$-elliptic.
 Thus, we may use the conclusion of  the previous 
paragraph for $K=\triangle_i$ obtaining that $\triangle_i$ is self-adjoint parametrix possessing, and consequently, that $D^{\bullet}$
 is a self-adjoint parametrix possessing complex. We may utilize the theorems derived in this paper.
Namely, using Theorems 5 and 6 and Corollary 7, one derives the items 2, 3 and 4.
According to Theorem 11 in \cite{KryslHodge}, $H^i(D^{\bullet},A)$ is a finitely generated $A$-module and a Banach space with respect to the 
quotient norm $| \,|_q.$ Thus due to Remark 1, $H^i(D^{\bullet},A)$ equipped with the canonical quotient structure,  is a finitely generated Hilbert $A$-module.
Using Theorem 5.9 in Frank, Larson \cite{FL}, the cohomology $H^i(D^{\bullet},A)$  is projective due to the unitality of  $A.$ 
$\Box$

{\bf Remark 7:} Let us notice that the decomposition in item 2 of the previous theorem is meant with respect to the $A$-product 
$(,)_{\Gamma}.$ Also the adjoints are considered with respect to this $A$-product.

\end{document}